\documentclass[12pt,a4paper]{article}
\usepackage{latexsym}
\newtheorem{lem}{Lemma}
\newtheorem{cor}{Corollary}
\newtheorem{thm}{Theorem}
\newenvironment{pf}{\textbf{Proof\ }}{\hfill$\Box$\smallskip}
\newcommand{\ov}{\overline}
\newcommand{\co}{\mathcal{O}}
\newcommand{\Q}{\mathbf{Q}}
\newcommand{\Z}{\mathbf{Z}}
\newcommand{\al}{\alpha}
\newcommand{\be}{\beta}

\newcommand{\la}{\lambda}
\newcommand{\La}{\Lambda}
\newcommand{\om}{\omega}
\title{Cubic identities for theta series in three variables}
\author{Robin Chapman\\
School of Mathematical Sciences\\ University of Exeter\\
Exeter, EX4 4QE, UK\\ \texttt{rjc@maths.ex.ac.uk}}
\date{20 September 2000}
\begin{document}
\maketitle

\section{Introduction}

In \cite{BB} (see also \cite{BBG}) Borwein and Borwein proved the
identity
\begin{equation}\label{originalid}
a(q)^3=b(q)^3+c(q)^3
\end{equation}
where
$$a(q)=\sum_{m,n\in\Z} q^{m^2+mn+n^2},$$
$$b(q)=\sum_{m,n\in\Z} \om^{m-n}q^{m^2+mn+n^2}$$
and
$$c(q)=\sum_{m,n\in\Z} q^{(m+1/3)^2+(m+1/3)(n+1/3)+(n+1/3)^2}$$
where $\om=\exp(2\pi i/3)$. We call these functions theta series for
convenience. Subsequently
Hirschhorn, Garvan and J.~Borwein \cite{HGB} proved the corresponding
identity for two-variable analogues of these theta series.
Sol\'e \cite{S} (see also \cite{SL}) gave a new proof
of (\ref{originalid}) using a lattice having the structure of a
$\Z[\om]$-module. Here we introduce three-variable analogues
of the theta series $a(q)$, $b(q)$ and $c(q)$, and adapt Sol\'e's
method to prove corresponding identities for them.

\section{Theta series}

We introduce our three-variable theta series as sums over elements
of the Eisenstein field $\Q(\sqrt{-3})$.

Let $K=\Q(\sqrt{-3})$ and let $\co=\Z[\om]$ be its ring of integers,
where $\om=\frac12(-1+\sqrt{-3})=\exp(2\pi i/3)$. Write
$\la=\om-\om^2=\sqrt{-3}$. For $\al\in K$ define $T(\al)=\al+\ov\al$,
the trace of $\al$. The element $\la$ generates a prime ideal of $\co$
of norm~3; the inclusion $\Z\to\co$ induces an isomorphism
$\Z/3\Z\cong \co/\la\co$.
Hence we can unambiguously define,
for $\al\in\co$, $\chi(\al)=\om^a$ where $a\in\Z$
and $\al\equiv a$ (mod~$\la\co$).

We now define our theta series. We start with
$$a(q,z,w)=\sum_{\al\in\co}q^{|\al|^2}z^{T(\al)}w^{T(\al/\la)}.$$
Next for any integer $k$ define
$$b_k(q,z,w)=\sum_{\al\in\co}\chi(\al)^kq^{|\al|^2}z^{T(\al)}w^{T(\al/\la)}.$$
It is apparent that $b_k(q,z,w)$ depends only on the congruence class of
$k$ modulo~3 and that $b_0(q,z,w)=a(q,z,w)$. We also define
$$c_k(q,z,w)=\sum_{\al\in\co+k/\la}q^{|\al|^2}z^{T(\al)}w^{T(\al/\la)}.$$
Again $c_k(q,z,w)$ depends only on the congruence class of
$k$ modulo~3 and $c_0(q,z,w)=a(q,z,w)$.

We observe some symmetry properties of these functions.
\begin{lem}
We have
\begin{equation}\label{aSymm}
a(q,z,w)=a(q,z,w^{-1})=a(q,z^{-1},w^{-1})=a(q,z^{-1},w),
\end{equation}
\begin{equation}\label{bSymm}
b_k(q,z,w)=b_k(q,z,w^{-1})=b_{-k}(q,z^{-1},w^{-1})=b_{-k}(q,z^{-1},w)
\end{equation}
and
\begin{equation}\label{cSymm}
c_k(q,z,w)=c_{-k}(q,z,w^{-1})=c_{-k}(q,z^{-1},w^{-1})=c_k(q,z^{-1},w).
\end{equation}
\end{lem}
\begin{pf}
We replace $\al$ in the definition of each series in turn by $\ov\al$,
$-\al$ and $-\ov\al$. It helps to note that $T(\ov\al)=T(\al)$,
$T(\ov\al/\la)=-T(\al/\la)$, $\chi(\ov\al)=\chi(\al)$,
$\chi(-\al)=\chi(\al)^{-1}$ and $\ov{k/\la}=-k/\la$.
Of course (\ref{aSymm}) is a special case of both
(\ref{bSymm}) and (\ref{cSymm}).
\end{pf}

From (\ref{bSymm}) and (\ref{cSymm}) we see that
$b_1(q,1,1)=b_{-1}(q,1,1)$ and $c_1(q,1,1)=c_{-1}(q,1,1)$. We write
$$a(q)=a(q,1,1),\quad b(q)=b_1(q,1,1)\quad\mathrm{and}\quad c(q)=c_1(q,1,1).$$
We shall soon see that this agrees with our previous definition.

We show that these functions specialize to the two-variable functions
introduced in \cite{HGB}. First of all, each element $\al\in\co$ can be
uniquely written as $\al=n\om-m\om^2$. Then $T(\al)=m-n$ and
$$|\al|^2=(n\om-m\om^2)(n\om^2-m\om)=m^2+mn+n$$
and so
$$a(q,z,1)=\sum_{m,n\in\Z}q^{m^2+mn+n^2}z^{m-n}$$
which is denoted as $a(q,z)$ in \cite{HGB}. In particular
$$a(q,1,1)=\sum_{m,n\in\Z}q^{m^2+mn+n^2}$$
in agreement with the original definition.
Also $|{-\om\al}|^2=|\al|^2$
and $T(-\om\al)=T((m-n\om^2)/\la)=n$. Hence
$$a(q,1,z)=\sum_{\al\in\co}q^{|{-\om\al}|^2}z^{T(-\om\al/\la)}
=\sum_{m,n\in\Z}q^{m^2+mn+n^2}z^n$$
which is denoted as $a'(q,z)$ in \cite{HGB}. Now $\chi(-\om\al)=\om^{m-n}$
and similarly
$$b_1(q,1,z)=\sum_{m,n\in\Z}\om^{m-n}q^{m^2+mn+n^2}z^n$$
which is denoted as $b(q,z)$ in \cite{HGB}. In particular
$$b_1(q,1,1)=\sum_{m,n\in\Z}\om^{m-n}q^{m^2+mn+n^2}.$$
Note that
$b_1(q,1,z)=b_{-1}(q,1,z)$ by (\ref{bSymm}).
Finally $\frac{1}{3}(\om-\om^2)=-1/\la$ and so
$$c_{-1}(q,z,1)=\sum_{m,n\in\Z}q^{(m+1/3)^2+(m+1/3)(n+1/3)+(n+1/3)^2}
z^{m-n}$$
and this is denoted by $q^{1/3}c(q,z)$ in \cite{HGB}. 
Again note that
$c_{-1}(q,z,1)=c_1(q,z,1)$ by (\ref{bSymm}).
In particular
$$c_1(q,1,1)=c_{-1}(q,1,1)=
\sum_{m,n\in\Z}q^{(m+1/3)^2+(m+1/3)(n+1/3)+(n+1/3)^2}.$$

\section{Identities}

Our main result is a generalization of (1.25) in \cite{HGB}.

\begin{thm}\label{thm1}
For each integer~$k$,
\begin{eqnarray}\label{cCubed}
3c_k(q,z,w)^3&=&a(q,w,z^{-3})a(q)^2\nonumber\\
&&{}+\om^k b_1(q,w,z^{-3})b(q)^2+\om^{-k}b_{-1}(q,w,z^{-3})b(q)^2\nonumber\\
&&{}+c_1(q,w,z^{-3})c(q)^2+c_{-1}(q,w,z^{-3})c(q)^2.
\end{eqnarray}
In particular
\begin{eqnarray}\label{aCubed}
3a(q,z,w)^3&=&a(q,w,z^{-3})a(q)^2
+b_1(q,w,z^{-3})b(q)^2+b_{-1}(q,w,z^{-3})b(q)^2\nonumber\\
&&{}+c_1(q,w,z^{-3})c(q)^2+c_{-1}(q,w,z^{-3})c(q)^2.
\end{eqnarray}
\end{thm}
\begin{pf}
Cubing the definition of $c_k(q,z,w)$ gives
\begin{eqnarray}\label{cCube}
c_k(q,z,w)^3=\sum_{\al_0,\al_1,\al_2\in\co+k/\la}
q^{|\al_0|^2+|\al_1|^2+|\al_2|^2}
z^{T(\al_0+\al_1+\al_2)}w^{T((\al_0+\al_1+\al_2)/\la)}.
\end{eqnarray}
This is a sum over triples $\al=(\al_0,\al_1,\al_2)$ where
$\al$ runs through a certain subset of
$$\La=\co^3+\Z(1/\la,1/\la,1/\la).$$

We partition the group $\La$ into various cosets. If
$\al\in\La$ then $\al_0+\al_1+\al_2\in\co$. For integers $j$
and $k$ let
$$\La_{j,k}=\{\al\in\co^3+k(1/\la,1/\la,1/\la):
\al_0+\al_1+\al_2\equiv j\pmod{\la}\}.$$
Then $\La_{j,k}$ depends only on the integers $j$ and $k$ modulo 3.
They are the nine cosets of the subgroup $\La_{0,0}$ of $\La$.
Define, for $\al\in K^3$,
$$|\al|^2=|\al_0|^2+|\al_1|^2+|\al_2|^2$$
and
$$\Phi(\al)=q^{|\al|^2}
z^{T(\al_0+\al_1+\al_2)}w^{T((\al_0+\al_1+\al_2)/\la)}.$$
Then
\begin{equation}\label{cCube2}
c_k(q,z,w)^3=\sum_{\al\in\La_{0,k}\cup\La_{1,k}\cup\La_{-1,k}}
\Phi(\al).
\end{equation}

We now consider the matrix
$$M=\frac{1}{\la}\left(\begin{array}{ccc}
1&1&1\\
1&\om&\om^2\\
1&\om^2&\om
\end{array}\right).$$
It is straightforward to check that $\La_{j,k}M=\La_{-k,j}$.
Also $M$ is a unitary matrix so that if $\be=\al M$ then
$|\be|^2=|\al|^2$.
Thus
$$\Phi(\al)=q^{|\be|^2}z^{T(\la\be_0)}w^{T(\be_0)}=
q^{|\be|^2}z^{-3T(\be_0/\la)}w^{T(\be_0)}.$$
From (\ref{cCube2}) we get
\begin{equation}\label{cCube3}
c_k(q,z,w)^3=\sum_{\be\in\La_{-k,0}\cup\La_{-k,1}\cup\La_{-k,-1}}
q^{|\be|^2}z^{-3T(\be_0/\la)}w^{T(\be_0)}.
\end{equation}
We split this sum into sums over each of the three cosets $\La_{-k,j}$.

Consider $\La_{-k,0}$. This can be written as
$$\La_{-k,0}=\{\be\in\co^3:\chi(\be_0)\chi(\be_1)\chi(\be_2)=\om^{-k}\}.$$
Hence
\begin{eqnarray}\label{firstcoset}
&&3\sum_{\be\in\La_{-k,0}}q^{|\be|^2}z^{T(-3\be_0/\la)}w^{T(\be_0)}\nonumber\\
&=&\sum_{\be\in\co^3}q^{|\be|^2}z^{T(-3\be_0/\la)}w^{T(\be_0)}\nonumber\\
&&{}+\om^k\sum_{\be\in\co^3}\chi(\be_0)\chi(\be_1)\chi(\be_2)
q^{|\be|^2}z^{T(-3\be_0/\la)}w^{T(\be_0)}\nonumber\\
&&{}+\om^{-k}\sum_{\be\in\co^3}\chi(\be_0)^{-1}\chi(\be_1)^{-1}\chi(\be_2)^{-1}
q^{|\be|^2}z^{T(-3\be_0/\la)}w^{T(\be_0)}\nonumber\\
&=&a(q,w,z^{-3})a(q)^2+\om^kb_1(q,w,z^{-3})b(q)^2
+\om^{-k}b_{-1}(q,w,z^{-3})b(q)^2.
\end{eqnarray}

To aid with the remaining cosets consider the matrix
$$N=\left(\begin{array}{ccc}
1&0&0\\
0&\om&0\\
0&0&\om
\end{array}\right).$$
Then $N$ is unitary and one may easily check that $\La_{j,k}N=\La_{j+k,k}$.
As $N$ does not alter the first coordinate of a triple $\be\in\ K^3$ then
for $k=\pm1$
$$\sum_{\be\in\La_{j,k}}q^{|\be|^2}z^{-3T(\be_0/\la)}w^{T(\be_0)}$$
is independent of~$j$. Hence for $k=\pm1$
\begin{eqnarray}\label{restcosets}
3\sum_{\be\in\La_{j,k}}q^{|\be|^2}z^{-3T(\be_0/\la)}w^{T(\be_0)}
&=&\sum_{\be\in(\co+k/\la)^3}q^{|\be|^2}z^{-3T(\be_0/\la)}w^{T(\be_0)}
\nonumber\\
&=&c_k(q,w,z^{-3})c(q)^2.
\end{eqnarray}
From (\ref{cCube3}), (\ref{firstcoset}) and (\ref{restcosets})
we obtain (\ref{cCubed}). The $k=0$ case of (\ref{cCubed}) is (\ref{aCubed}).
\end{pf}

\begin{cor}\label{cor1}
We have
\begin{eqnarray}\label{mainconclusion}
2a(q,z,w)^3&=&b_1(q,w,z^{-3})b(q)^2+b_{-1}(q,w,z^{-3})b(q)^2\nonumber\\
&&{}+c_1(q,z,w)^3+c_2(q,z,w)^3.
\end{eqnarray}
Also
\begin{equation}\label{qconclusion}
a(q)^3=b(q)^3+c(q)^3.
\end{equation}
\end{cor}
\begin{pf}
To obtain (\ref{mainconclusion}) subtract the sum the $k=1$ and $k=-1$
cases of (\ref{cCubed}) from twice (\ref{aCubed}).
To obtain (\ref{qconclusion}), either substitute $z=w=1$ in
(\ref{mainconclusion}), or make this substitution in either (\ref{aCubed})
or (\ref{cCubed}).
\end{pf}

Another particular case is obtained by setting $w=1$ to give
$$a(q,z,1)^3=b_1(q,1,z^3)b(q)^2+c_1(q,z,1)^3$$
(using (\ref{cSymm})) which is (1.25) in \cite{HGB}.

A variant of the argument of Theorem~\ref{thm1} gives the following
result.

\begin{thm}\label{thm2}
For each $k$,
\begin{eqnarray}\label{cCubed2}
3c_k(q,z,w)c_k(q^2,z^2,w^2)
&=&a(q,w,z^{-3})a(q^2)\nonumber\\
&&{}+\om^k b_1(q,w,z^{-3})b(q^2)+\om^{-k}b_{-1}(q,w,z^{-3})b(q^2)\nonumber\\
&&{}+c_1(q,w,z^{-3})c(q^2)+c_{-1}(q,w,z^{-3})c(q^2).
\end{eqnarray}
In particular
\begin{eqnarray}\label{aCubed2}
3a(q,z,w)a(q^2,z^2,w^2)
&=&a(q,w,z^{-3})a(q^2)\nonumber\\
&&{}+b_1(q,w,z^{-3})b(q^2)+b_{-1}(q,w,z^{-3})b(q^2)\nonumber\\
&&{}+c_1(q,w,z^{-3})c(q^2)+c_{-1}(q,w,z^{-3})c(q^2).
\end{eqnarray}
\end{thm}
\begin{pf}
As the proof follows closely the proof of Theorem~\ref{thm1},
we shall suppress most of the details.

Let $V=\{(\al_0,\al_1,\al_2)\in K^3:\al_1=\al_2\}$.
The space $V$ is stable under the action of the matrices
$M$ and $N$. The key is to rewrite the proof of Theorem~\ref{thm1}
restricting the summations to triples in~$V$. We start by noting that
$$c_k(q,z,w)c_k(q^2,z^2,w^2)
=\sum_{\al\in(\co+1/\la)^3\cap V}\Phi(\al).$$
This gives
$$c_k(q,z,w)c_k(q^2,z^2,w^2)
=\sum_{\be\in(\La_{-k,0}\cap V)\cup(\La_{-k,1}\cap V)\cup(\La_{-k,-1}\cap V)}
q^{|\be|^2}z^{-3T(\be_0/\la)}w^{T(\be_0)}.$$
We then get
\begin{eqnarray*}
&&3\sum_{\be\in\La_{-k,0}\cap V}q^{|\be|^2}z^{T(-3\be_0/\la)}w^{T(\be_0)}\\
&=&a(q,w,z^{-3})a(q^2)+\om^kb_1(q,w,z^{-3})b(q^2)
+\om^{-k}b_{-1}(q,w,z^{-3})b(q^2)
\end{eqnarray*}
and for $j=\pm1$
$$3\sum_{\be\in\La_{j,k}\cap V}q^{|\be|^2}z^{-3T(\be_0/\la)}w^{T(\be_0)}
=c_k(q,w,z^{-3})c(q^2).$$
The theorem then follows.
\end{pf}

\begin{cor}\label{cor2}
We have
\begin{eqnarray*}
2a(q,z,w)a(q^2,z^2,w^2)&=&b_1(q,w,z^{-3})b(q)^2+b_{-1}(q,w,z^{-3})b(q)^2\\
&&{}+c_1(q,z,w)^3+c_2(q,z,w)^3.
\end{eqnarray*}
Also
$$a(q)a(q^2)=b(q)b(q^2)+c(q)c(q^2).$$
\end{cor}
\begin{pf}
This follows from Theorem~\ref{thm2} in exactly the same
way that Corollary~\ref{cor1} follows from Theorem~\ref{thm1}.
\end{pf}

Another special case is
$$a(q,z,1)a(q^2,z^2,1)=b_1(q,1,z^3)b(q^2)+c_1(q,z,1)c_1(q^2,z^2,1)$$
which is (1.26) in \cite{HGB}.


\begin{thebibliography}{99}

\bibitem{BB}
J.~M.~Borwein \& P.~B.~Borwein, 
A cubic counterpart of Jacobi's identity and the AGM,
\emph{Trans.\ Amer.\ Math.\ Soc.} \textbf{323} (1991), 691--701.

\bibitem{BBG}
J.~M.~Borwein, P.~B.~Borwein \& F.~G.~Garvan,
Some cubic modular identities of Ramanujan,
\emph{Trans.\ Amer.\ Math.\ Soc.} \textbf{343} (1994), 35--47.

\bibitem{HGB}
M.~Hirschhorn, F.~Garvan \& J.~Borwein,
Cubic analogues of the Jacobian theta function $\theta(z,q)$,
\emph{Canad.\ J.\ Math.} \textbf{45} (1993), 673--694.

\bibitem{S}
P.\ Sol\'e, 
$D_4$, $E_6$, $E_8$ and the AGM,
\emph{Springer Lecture Notes in Computer Science}
\textbf{948} (1995), 448--455.

\bibitem{SL}
P.\ Sol\'e \& P.\ Loyer,
$U_n$ lattices, construction $B$, and AGM iterations,
\emph{European J.\ Combin.} \textbf{19} (1998), 227-236.

\end{thebibliography}
\end{document}